# Existence of boundary flex control for the systems governed by Boussinesq equation with the press boundary condition and mixed boundary condition

Gol Kim [a]

[a] Center of Natural Sciences, University of Sciences, DPR Korea

**Abstract.** In this paper, the boundary flex control problem of non stationary equation governing the coupled mass and heat flow of a viscous incompressible fluid in a generalized Boussinesq approximation by assuming that viscosity and heat conductivity are dependent on temperature has been studied. The boundary condition for velocity of fluid is non -standard boundary condition: specifically the case where dynamical pressure is given on some part of the boundary and the boundary condition for temperature of fluid is mixed boundary condition has been considered.
First, we have proved the existence of existence of the weak solution for state equation.Then the optimal condition has been proved the existence of optimal control.

## 1. Introduction

The extremal problem for Navier-Stokes equation and optimal control of fluid dynamical equation are studied by various authors (for example; [6-10], [17-19]).

In [22] the existence of time optimal controls for the Boussinesq equation has been obtained and derived Pontryagin's maximum principle of time optimal control problem governed by the Boussinesq equation

In [23] an optimal control problem governed by a system of nonlinear partial differential equations modeling viscous incompressible flows submitted to variations of temperature has been consider. A generalized Boussinesq approximation has been used. The existence of the optimal control as well as first order optimality conditions of Pontryagin type by using the Dubovitskii-Milyutin formalism has been obtained.

In [24] the stationary Boussinesq equations describing the heat transfer in the viscous heat-conducting fluid under inhomogeneous Dirichlet boundary conditions for velocity and mixed boundary conditions for temperature are considered. The optimal control problems for these equations with tracking-type functionals are formulated. A local stability of the concrete control problem solutions with respect to some disturbances of both cost functionals and state equation is proved.

In [25] the boundary control problems of the model of heat and mass transfer in a viscous incompressible heat conducting fluid has been considered. The model consists of the Navier-Stokes equations and the convection-diffusion equations for the substance concentration and the temperature that are nonlinearly related via buoyancy in the Oberbeck–Boussinesq approximation and via convective mass and heat transfer.

In [25] control problems for stationary magnetohydrodynamic equations of a viscous heat-conducting fluid under mixed boundary conditions has been considered.

In [15, 16] the Karhunen-Loeve Galerkin method for the inverse problems of Boussinesq equation have been studied.

In [20] the problem of stabilization of the Boussinesq equation via internal feedback controls has been studied. In [21] the problem of local internal controllability of the Boussinesq system has been studied.

In this paper, the boundary flex optimal control for the evolution equation governing the coupled mass and heat flow of a viscous incompressible fluid in a generalized Boussinesq approximation by assuming that viscosity and heat conductivity are dependent on temperature has been studied. The boundary condition for velocity of fluid is non -standard boundary condition: specifically the case where dynamical pressure is given on some part of the boundary will be considered. The boundary condition for temperature of fluid is mixed boundary condition.

The existence of the optimal control has been proved. Then the optimal condition has been derived.

Let $\Omega \subset R^N$ (N=2, 3) be a bounded domain with smooth boundary $\Gamma$. Let $\Gamma$ be divided by into two parts $\Gamma_1, \Gamma_2$ such that $\Gamma = \Gamma_1 \cup \Gamma_2 (\Gamma_1 \cap \Gamma_2 = \varnothing, \Gamma_2 \neq \varnothing)$ T. ($0<T<\infty$) is given number.

We denote $Q = \Omega \times (0, T)$, $\Sigma_i = \Gamma_i \times (0, T)$ (i = 1,2), $\Sigma = \Gamma \times (0, T)$.

We assume that the state of control systems is given by non-stationary Boussinesq equation with the dynamical pressure condition and mixed boundary condition on some part of the boundary as follows:

$$\begin{cases} \dfrac{\partial z}{\partial t} - \nu \Delta z + (z, \nabla)z + \beta g\, w = -\operatorname{grad} \pi; & Q \\ \operatorname{div} z = 0; & Q \\ \dfrac{\partial w}{\partial t} - k \Delta \nabla w + (z, \nabla)w = 0; & Q \\ z_\varsigma = 0,\ \pi + \dfrac{1}{2}|z|^2 = v_1,\ w = 0; & \Sigma_1 \\ z = 0,\ -k(w)\dfrac{\partial w}{\partial n} = v_2; & \Sigma_2 \\ z(0) = z_0,\ w(0) = w_0; & \Omega \end{cases} \quad (1.1)$$

where $\Omega \subset R^N$ (N=2, 3) is a bounded domain with smooth boundary $\Gamma$. $\Gamma$ is divided by into two parts $\Gamma_1, \Gamma_2$ such that $\Gamma = \Gamma_1 \cup \Gamma_2 (\Gamma_1 \cap \Gamma_2 \neq \Phi)$, T ($0<T<\infty$) is given number.

We denote $Q = \Omega \times (0, T)$, $\Sigma_i = \Gamma_i \times (0, T)$ (i = 1,2), $\Sigma = \Gamma \times (0, T)$ and $n$ note the outer normal vector to $\Gamma$. In the Eqs. (1)-(6) $z(x, t) \in R^N$ denotes the velocity of the fluid at point $x \in \Omega$ at time $t \in [0, T]$; $\pi(x, t) \in R$ is the hydrostatic pressure; $w(x, t) \in R$ is temperature; g is the gravitational vector, and $\nu > 0$ and $k > 0$ are kinematic viscosity and thermal conductivity, respectively; $\beta$ is a positive constant associated to the coefficient of volume expansion; $v_1$ and $v_2$ are the given functions on $\Sigma_1$ and $\Sigma_2$ respectively. In Eqe.(1) $\beta > 0$ is the coefficient of volume expansion and g is the gravitational function.

The expressions $\nabla, \Delta$ and div denote the gradient, Laplacian and divergence operators, respectively (sometimes, we will also denote the gradient operator by grad); $i$ th component in Cartesian coordinates of $(z, \nabla)z$ is given by

$$((z, \nabla)z)_i = \sum_{j=1}^{N} z_j \frac{\partial z_i}{\partial x_j}, \text{ also } (z, \nabla)w = \sum_{j=1}^{N} z_j \frac{\partial w}{\partial x_j}$$

In the boundary condition (4) $z_\varsigma = z - z_n$, $z_n = (z \cdot n)n$.

There are the results of research of Boussinesq equation and the generalized Boussinesq system with nonlinear thermal diffusion in [1-4, 23]. But boundary conditions of those papers are homogeneous.

We assume that the cost functional J[v] is given as following:

$$J[v, y] = N_1 \int_0^T \int_{\Gamma_1} r_1(s,t) z_n(s,t)\, ds\, dt + N_2 \int_{\Gamma_2} r_2(s,t) \frac{\partial w}{\partial n}\, ds\, dt$$

where $N_1, N_2 > 0$ are given real numbers and $r_1(x, t) \in [L^2(\Sigma_1)]^N$, $r_2(x, t) \in L^2(\Sigma_2)$ are given functions. And

$\int_{\Gamma_1} r_1(s, t) z_n(s, t)\, ds = \int_{\Gamma_1} (r_1(s, t), z_n(s, t))\, ds$ and $(r_1(s, t), z_n(s, t))$ denote the scalar product in $[L^2(\Sigma_1)]^N$.

Then, the problem that we are going to considered is to find the $v_* \in U_a$ satisfying:

$$\inf_{v \in U_a} J[v] = J[v_*] \qquad (1.2)$$

Here, we assume that the admissible control sets $U_\alpha = U_{1\alpha} \times U_{2\alpha}$ are defined such as:

$$U_{1\alpha} = \{v_1 \mid v_1 \in L^2(0,T;(L^2(\Gamma_1))^N), 0 < \alpha_1(x,t) \le v_1(x,t) \le \beta_1 (almost\ everywhere)\} \qquad (1.3)$$

$$U_{2\alpha} = \{v_2 \mid v_2 \in L^2(0,T;L^2((\Gamma_2))), 0 < \alpha_2 \le v_2(x,t) \le \beta_2 (almost\ everywhere)\} \qquad (1.4)$$

$\alpha_i(x,t), \beta_i(x,t) (i=1,2)$ are given functions in the function space $L^2(\Sigma_i)$.
$v_1(x,t) = (v_{11}(x,t), v_{12}(x,t), \cdots, v_{1N}(x,t))$ and expression $0 < \alpha_1(x,t) \le v_1(x,t) \le \beta_1$ means that
$0 < \alpha_1(x,t) \le v_{1i}(x,t) \le \beta_1\ (\forall i \in N)$

We denote $v = \{v_1, v_2\}$ and $y(v) = \{z(v), w(v)\}$.

The established optimization problem (1.1)-(1.4) is an optimal boundary flex control problem.

For the convenient, we have assumed that control parameters are the flex pressure $v_1$ on the boundary $\Sigma_1$ and the heat flex $v_2$ on the boundary $\Sigma_2$

To illustrate the example of extremal condition (1.2), we can take functions $r_1(x,t)$ and $r_2(x,t)$ as following;

$$r_1(x,t) = \begin{cases} I & ;x \in \overline{\Gamma}_1 \subset \Gamma_1 \\ 0 & ;x \notin \overline{\Gamma}_1 \subset \Gamma_1 \end{cases}, \qquad r_2(x,t) = \begin{cases} 1 & ;x \in \overline{\Gamma}_2 \subset \Gamma_2 \\ 0 & ;x \notin \overline{\Gamma}_2 \subset \Gamma_2 \end{cases}$$

Where, I is a unit vector. Then, the optimization problem (1-2) is described the problem that fluid flex passed the boundary $\overline{\Gamma}_1$ under the restriction for the flex pressure and heat flex passed the boundary $\overline{\Gamma}_2$ under the restriction for the heat flex must minimize.

## 2. Preliminaries

In this article the functions are either R or $R^N$ ($N=2$ or $N=3$) and as usual simplification, sometimes we will not distinguish them in our notations; the difference will be clear from the context. The $L^2(\Omega)$-product and norm are denoted by $(\cdot,\cdot)$ and $|\cdot|$ respectively: the $H^m(\Omega)$ norm is denoted by $\|\cdot\|_m$. Here $H^m(\Omega) = W^{m,2}(\Omega)$ is the usual Sobolev spaces (see [1] for their properties); $H^{-1}(\Omega)$ denotes the dual spaces of $H_0^1(\Omega)$. $H^0(\Omega)$ is the same as $L^2(\Omega)$ and $\|\cdot\|_0$ is the same as $L^2$-norm $|\cdot|$. $D(0,T)$ is the class of $C^\infty$-functions with compact support in $(0,T)$. $D(0,T)'$ are its associated spaces of distribution.

If B is a Banach space, we denoted by $L^q(0,T;B)$ the Banish space of the B-valued functions defined in the interval $(0,T)$ that are $L^q$-integrable.

Now we introduce some spaces such as;

$$\begin{cases} D = \{\phi \mid \phi \in (c^\infty(\Omega))^N, div\phi(x) = 0 (x \in \Omega), \phi_\varsigma(x) = 0 (x \in \Gamma_2)\} \\ H = \text{completion of D under the } [L^2(\Omega)]^N\text{-norm} \\ V = \text{completion of D under the } [H^1(\Omega)]^N - \text{norm} \\ D_{\Gamma_1} = \{\phi : \phi \in C^\infty(\Omega), \phi(x) = 0 (x \in \Gamma_1)\} \\ \widetilde{H} = \text{closure of } D_{\Gamma_1} \text{ in } L^2(\Omega) \\ W = \text{closure of } D_{\Gamma_1} \text{ in } H^1(\Omega) \end{cases} \qquad (2.1)$$

Naturally, the norm of H or $\widetilde{H}$ is also denoted by $|\cdot|$, and the norm of $V$ or $W$ is denoted by $\|\cdot\|$ as well. The dual product between $V^*$ and $V$ or $W^*$ and $W$ (also the inner product in $H^{-1}(\Omega)$ and $H_0^1(\Omega)$) are denoted by $< \cdot, \cdot >$.

Now, generally, we assume that $v_1 \in L^2(0,T;(H^{-1/2}(\Gamma_1))^N)$, $v_2 \in L^2(0,T;H^{-1/2}(\Gamma_2))$.
Then, we shall prove the existence of the weak solution for state equation (1.1)
Suppose that $\{z,w\}$ is a classical solution of (1.1). Multiplier the first equation of (1.1) by $\psi \in V$, integrate by parts over $\Omega$ and take the boundary condition into account to get

$$\frac{d}{dt}(z,\psi) + v\, a_1(z,\psi) + b(z,z,\psi) + (\beta g\, w,\psi) = <v_1,\psi_n>_{\Gamma_1}, \quad \psi_n = (\psi \cdot n)n. \tag{2.2}$$

Multiplier the second equation of (1.1) by $\varphi \in W$, integrate by parts over $\Omega$, and take the boundary conditions into account to obtain

$$\frac{d}{dt}(w,\varphi) + k a_2(w,\varphi) + c(z,w,\varphi) = <v_2,\varphi>. \tag{2.3}$$

Now let $\chi \in C^1[0,T]$ be a function such that $\chi(T) = 0$. Multiplier Equations (2.2) and (2.3) by $\chi$ respectively, and integrate by parts to yield

$$\int_0^T [-(z(t),\psi\chi'(t))dt + \int_0^T [va_1(z(t),\psi\chi(t)) + b(z(t),z(t),\psi\chi(t)) + \beta(w(t)g,\psi\chi(t))]dt$$

$$= \int_0^T \int_{\Gamma_1} (v_1(t),\psi_n\chi(t))dsdt + (z_0,\psi)\chi(0) \tag{2.4}$$

$$\int_0^T [-(w(t),\varphi\chi'(t))dt + \int_0^T [ka_2(w(t),\varphi\chi(t)) + c(z(t),w(t),\varphi\chi(t))]dt$$

$$= \int_0^T \int_{\Gamma_{21}} (v_2(t),\varphi\chi(t))dsdt + (w_0,\varphi)\chi(0) \tag{2.5}$$

where $z(\cdot,t) = z(t)$, $w(\cdot,t) = w(t)$, $v_i(\cdot,t) = v_i(t)$, $i = 1,2$ by abuse of notation without confusion from the context, and

$$(z,\psi) = \sum_{j=1}^N \int_\Omega z_j(x,\cdot)\psi_j(x)dx, \quad (w,\varphi) = \int_\Omega w(x,\cdot)\varphi(x)dx$$

$$a_1(u,\psi) = (rotz(x,\cdot),rot\psi) = \int_\Omega (\nabla \times u)\cdot(\nabla \times \phi)dx \quad b(z,z,\phi) = \int_\Omega (rotz(x,\cdot)\times z(x,\cdot))\psi(x)dx,$$

$$a_2(w,\varphi) = \sum_{j=1}^N \int_\Omega \frac{\partial w(x,\cdot)}{\partial x_j} \cdot \frac{\partial \varphi}{\partial x_j}dx, \quad c(z,w,\varphi) = \sum_{j=1}^N \int_\Omega z_j(x,\cdot)\frac{\partial w(x,\cdot)}{\partial x_j}\varphi(x)dx$$

The following Lemmas 2.1-2.3 can be obtained by the Sobolev inequalities and the compactness theorem. We can also refer to theorem 1.1 of [31] on page 107 and lemmas 1.2, 1.3 of [31] on page 109 (see also chapter 2 of [13]). The similar arguments can be also found in lemmas 1, 5 of [30].

**Lemma 2.1.** The bilinear forms $a_1(\cdot,\cdot)$ and $a_2(\cdot,\cdot)$ are coercive over V and W respectively. That is, there exist constants $c_1, c_1' > 0$ such that

$$a_1(z,z) \geq c_1 \|z\|^2, \forall z \in V \text{ and } a_2(w,w) \geq c_1' \|z\|^2, \forall \in W$$

**Lemma 2.2.** The trilinear forms $b(\cdot,\cdot,\cdot)$ is a linear continuous functional on $[H^1(\Omega)]^N$. That is, there exist a constant c2 > 0 such that

$$|b(u,v,w)| \leq c_2 \|u\| \cdot \|v\| \cdot \|w\|, \qquad \forall u,v,w \in [H^1(\Omega)]^N$$

Moreover, the following properties hold true

( i ) b(u, v, v) = 0, $\forall$u, v $\in$ V
(ii) b(u, v, w) = $-$b(u, w, v), $\quad \forall$u $\in$ V, $\forall$v, w $\in [H^1(\Omega)]^N$

(iii) If $u_m \to u$ weakly on V and $v_m \to v$ strongly on H, then
$$b(u_m, v_m, w) \to b(u, v, w), \forall u, v \in V \quad \forall w \in V$$

**Lemma 2.3.** The tri-linear form $c(\cdot,\cdot,\cdot)$ is a linear continuous functional defined on $V \times W \times W$. That is, there exist a constant $c_3 > 0$ as following:
$$|c(z, w, \varphi)| \leq c_3 \|z\| \cdot \|w\| \cdot \|\varphi\|, \quad \forall z \in V, \quad \forall w, \varphi \in W$$

Moreover, the following properties hold true

(i) $c(z, w, w) = 0$, $\forall z \in V$, $\forall w \in W$

(ii) $c(z, w, \varphi) = -c(z, \varphi, w)$, $\forall z \in V, \forall w, \varphi \in W$

(iii). If $z_m \to z$ weakly on V and $w_m \to w$ strongly on $\tilde{H}$, then
$$c(z_m, w_m, \varphi) \to c(z, w, \varphi), \forall z \in V, w \in \tilde{H}, \varphi \in W.$$

**Definition 1.** Let $Y \equiv Z \times W = (L^2(0, T : V) \cap L^\infty(0, T : H)) \times (L^2(0, T : W) \cap L^\infty(0, T : \tilde{H}))$. Suppose that $v_1 \in (L^2(0, T; (H^{-1/2}(\Gamma_1))^N)$, $v_2 \in (L^2(0, T; H^{-1/2}(\Gamma_2)))$, $z_0 \in H$, $w_0 \in \tilde{H}$, $g \in L^\infty(\Omega)$. The pair $y = \{z, w\}$ is said to be a weak solution of (1.1) if it satisfies

$$\begin{cases} y = \{z, w\} \in Y, z' \in L^1(0, T, V^*), w' \in L^1(0, T : W^*) \\ (z', \psi) + \nu a_1(z, \psi) + b(z, z, \psi) + (\beta g w, \psi) = <v_1, \psi_n>_{\Gamma_1}, \forall \psi \in V \\ (w', \varphi) + k a_2(w, \varphi) + c(z, w, \varphi) = <v_2, \varphi>_{\Gamma_2}, \forall \varphi \in W \\ z(0) = z_0, w(0) = w_0 \end{cases} \quad (2.7)$$

Next, we reformulate Equation (2.7) into the operator equation. To this purpose, it is noticed that for a _fixed $\psi \in V$, the functional $\psi(\in V) \to a_1(z, \psi)$ is linear continuous. So there exists an $A_1 z \in V^*$ such that
$$<A_1 z, \psi> = a_1(z, \psi), \quad \forall \psi \in V \quad (2.8)$$

Similarly, for fixed $u, v \in V$, $w \in V \to b(u, v, w)$ is a linear continuous functional on V. Hence there exist a $B(u, v) \in V^*$ such that
$$<B(u, v), w> = b(u, v, w), \quad \forall w \in V \quad (2.9)$$

We denote $B(u) = B(u, u)$. Define
$$L_1(\psi) = (v_1, \psi_n)_{\Gamma_1} = \int_{\Gamma_1} v_1 \psi_n ds, \quad \forall \psi \in V. \quad (2.10)$$

Then, fixed $v_1 \in (L^2(0, T; (H^{-1/2}(\Gamma_1))^N)$, the functional $\psi(\in V) \to L_1(\psi) = (v_1, \psi_n)_{\Gamma_1}$ is linear continuous. So that there exist constant $c_4 > 0$ such that
$$\|L_1 \psi\| < c_4 \|\psi\|_V, \forall \psi \in V.$$

So there exists an $H_1 v_1 \in V^*$ such that
$$<H_1 v_1, \psi> = (v_1, \psi_n)_{\Gamma_1}, \forall \psi \in V \quad (2.11)$$

With these operators at hand, we can write the second equation of (2.7) as
$$\frac{dz}{dt} + \nu A_1 z + B(z) + \beta g w = H_1 v_1 \quad (2.12)$$

Similarly, we have
$$<A_2 w, \varphi> = a_2(w, \varphi), \quad <C(z, w), \varphi> = c(z, w, \varphi), \quad A_2 w, C(z, w) \in W^* \quad (2.13)$$

Define
$$L_2(\varphi) = (v_2, \varphi)_{\Gamma_2} = \int_{\Gamma_2} v_2 \varphi ds, v_2 \in L^2(\Gamma_2), \forall \varphi \in W \quad (2.14)$$

Then, fixed $v_2 \in (L^2(0, T; H^{-1/2}(\Gamma_2))$ the functional $\varphi(\in W) \to L_2(\varphi) = (v_2, \varphi)_{\Gamma_{21}}$ is linear continuous.

Then, the operator $L_2$ is a linear continuous functional defined on $W$ and so there exists constant $c_5 > 0$ such that
$$\|L_2\varphi\| \leq c_5 \|\varphi\|_W, \forall \varphi \in W.$$
So there exists an $H_2 v_2 \in W^*$ such that
$$<H_2 v_2, \varphi> = (v_2, \varphi)_{\Gamma_2}, \forall \varphi \in V \tag{2.15}$$
By these operators defined above, we can write the first equation of (2.7) as
$$\frac{dw}{dt} + kA_2 w + C(z, w) = H_2 v_2 \tag{2.16}$$
Combining (2.12) and (2.16), we can write (2.7) in the abstract evolution equation as follows:
$$\begin{cases} \frac{dz}{dt} + \nu A_1 z + B(z) + \beta gw = H_1 v_1, \\ \frac{dw}{dt} + kA_2 w + C(z, w) = H v_2, \\ z(0) = z_0, w(0) = w_0 \end{cases} \tag{2.17}$$

**Lemma 2.4.** If $z \in L^2(0,T;V)$, then $B(z) \in L^1(0,T;V^*)$; and if $w \in L^2(0,T;W)$, then $C(z,w) \in L^1(0,T;W^*)$.

Proof. By applying Höler inequality and compactness of embedding $H^1(\Omega) \subset L_4(\Omega)$, we obtain
$$|(B(z),\psi)| = |b(z,z,\psi)| \leq c_6' \|z\| \cdot \|z\|_{L_4} \cdot \| \|_{L_4} \leq c_6 \|z\|^2 \|\psi\|,$$
for some constants $c_6', c_6 > 0$ and hence $\|B(z)\|_{V^*} \leq c_6 \|z\|^2$ which shows that $B(z) \in L^1(0,T;V^*)$. The proof is complete.
Similarly, we have
$$|(C(z,w),\varphi)| = |-c(z,\varphi,w)| \leq c_7 \|z\|_{L_4} \cdot \|\varphi\| \cdot \|w\|_{L_4} \leq c_8 \|z\| \cdot \|\varphi\| \cdot \|w\| \leq c_9 (\|z\|^2 + \|w\|^2)\|\varphi\|$$
for some constants $c_8, c_9 > 0$ and hence $\|C(z,w)\| \leq c_9 (\|z\|^2 + \|w\|^2)$ which shows that $C(z,w) \in L^1(0,T;W^*)$ for all $\varphi \in W$. The proof is complete.

We specify the constants $c_i, i = 1,2,\cdots$ used in this section in the remaining part of the paper. The following Lemma 2.5 comes from theorem 2.2 of [12] on page 220

**Lemma 2.5.** Let $X_0, X, X_1$ be Hilbert spaces with the compact embedding relations
$$X_0 \subset X \subset X_1$$
Then, for arbitrary bounded set $K \subset R^1$, $\nu > 0$, embedding $H_K^\nu(R^1; X_0, X_1) \subset L^2(R^1; X)$ is compact, where
$$H_K^\nu(R^1; X_0, X_1) = \{v \in H_K^\nu(R^1; X_0, X_1) : \sup pv \subset K, \hat{p}_t^\nu v \in L^2(R; X_1)\}$$
$$H_K^\nu(R^1; X_0, X_1) = \{v \in L^2(R^1; X_0), \hat{p}_t^\nu v \in L^2(R^1; X_1)\}$$
$$\hat{p}_t^\nu v(\tau) = (2\pi i\tau)^\nu \hat{v}(\tau), \quad \hat{v}(\tau) = \int_{-\infty}^{+\infty} e^{-2\pi i\tau} v(t)dt$$
$$\|v\|^2_{H_K^\nu(R;X_0;X_1)} = \|v\|^2_{L^2(R^1;X_0)} + \||\tau|^2 \hat{v}\|^2_{L^2(R;X_1)}$$

## 3 Existence of the weak solution for state equation

This section discusses the existence of the weak solution defined by Definition 1 to Equation (1.1). The main idea is to construct a Galerkin approximation scheme and make some prior estimates.

Choose two orthogonal bases $\{u_j\}_{j=1}^\infty$ for $V$ and $\{\mu_j\}_{j=1}^\infty$ for W respectively. Construct the

Galerkin approximation solutions:
$$z_m(x,t) = \sum_{j=1}^{m} q_{jm}(t) u_j(x), \quad w_m(x,t) = \sum_{j=1}^{m} h_{jm}(t) \mu_j(x) \tag{3.1}$$

such that for all $j \in N^+$, $\{z_m, w_m\}$ satisfies

$$\begin{cases} (z'_m(t), u_j) + a_1(z_m(t), u_j) + b(z_m(t), z_m(t), u_j) + \beta(w_m(t)g, u_j) = (v_1(t), u_{jn})_{\Gamma_1} \\ (w'_m(t), \mu_j) + a_2(w_m(t), \mu_j) + c(z_m(t), w_m(t), \mu_j) = (v_2(t), \mu_j)_{\Gamma_{21}} \\ z_m(0) = z_{m0} \to z_0 \text{ in } H, w_m(0) = w_{m0} \to w_0 \text{ in } \widetilde{H}, \quad j = 1, 2, \cdots \end{cases} \tag{3.2}$$

where $z_{m0}$ is the orthogonal projection of $z_0$ in $H$ on the subspace spanned by $\{u_j\}_{j=1}^{\infty}$ and $w_{m0}$ is the orthogonal projection of $w_0$ in $\widetilde{H}$ on the subspace spanned by $\{\mu_j\}_{j=1}^{\infty}$.

Once again, we write $z_m(x,t) = z_m(\cdot, t)$, $w_m(x,t) = w_m(\cdot, t)$ by abuse of notation without the confusion from the context. It is seen that for any $m \in N^+$, system (3.2) is a system of nonlinear differential equations with the unknown variables $\{q_{jm}(t), h_{jm}(t)\}$ and the initial values $\{q_{jm}(0) = (z_0, u_j), h_{jm}(0) = (w_0, \mu_j)\}_{j=1}^{m}$. By the assumption, this initial value problem admits a solution in some interval $[0, t_m]$. We need to show that $t_m = T$ tm = T.

**Lemma 3.1.** Let $\{z_m, w_m\}$ be the sequence satisfying (3.2). Then there exists a subsequence of $\{z_m, w_m\}$, still denoted by itself without confusion, such that

$$z_m \to z \text{ weakly in } L^2(0,T;V) \text{ and } z_m \to z \text{ weakly star in } L^\infty(0,T;H), \tag{3.3}$$

where $z \in L^2(0,T;V) \cap L^\infty(0,T;H)$, and

$$w_m \to w \text{ weakly in } L^2(0,T;W) \text{ and } w_m \to w \text{ weakly star in } L^\infty(0,T;\widetilde{H}), \tag{3.4}$$

where $w \in L^2(0,T;W) \cap L^\infty(0,T;\widetilde{H})$,.

Proof. Sum for $j$ from 1 to m in (3.2) and apply Lemmas 2.2, 2.3, to get

$$\begin{cases} (z'_m(t), z_m(t)) + a_1(z_m(t), z_m(t)) + \beta(w_m(t)g, z_m) = (v_1(t), z_{mn})_{\Gamma_1} \\ (w'_m(t), w_m(t)) + a_2(w_m(t), w_m(t)) = (v_2(t), w_m)_{\Gamma_{21}} \end{cases} \tag{3.5}$$

By assumption (2.6), for any given $\varepsilon > 0$, we can get from (3.5) that

$$\frac{d}{dt}|z_m(t)|^2 + c_1 \|z_m(t)\|^2 \leq -\beta(w_m(t)g, z_m(t)) + (v_1(t), z_{mn})_{\Gamma_1} \leq \frac{1}{2}\beta\|g\|_\infty (\frac{1}{\varepsilon^2}\|w_m(t)\|^2 + \varepsilon^2\|w_m(t)\|^2)$$

$$+ \frac{1}{2}\varepsilon^2 \|z_{mn}(t)\|^2_{H^{1/2}(\Gamma_1)} + \frac{1}{2\varepsilon^2}\|v_1(t)\|^2_{H^{-1/2}(\Gamma_1)} \tag{3.6}$$

Here and in what follows, we denote $\|z_{mn}(t)\|^2_{(H^{1/2}(\Gamma_1))^N}$ and $\|v_1(t)\|^2_{(H^{-1/2}(\Gamma_1))^N}$ simply by $\|z_{mn}(t)\|^2_{H^{1/2}(\Gamma_1)}$ and $\|v_1(t)\|^2_{H^{-1/2}(\Gamma_1)}$ respectively by abuse of the notation.

By the trace theorem from $H^1(\Omega)$ to $H^{1/2}(\Gamma)$, there exists a constant $c_{10}$ such that

$$\|z_{mn}(t)\|_{H^{1/2}(\Gamma_1)} \leq c_{10} \|z_m(t)\|$$

Substitute above into (3.6) to yield

$$\frac{d}{dt}|z_m(t)|^2 + [c_1 - (\beta\|g\|_\infty \varepsilon^2 + c_{10}\varepsilon^2)]\|z_m(t)\|^2$$

$$\leq \beta\|g\|_\infty \frac{1}{\varepsilon^2}\|w_m(t)\|^2 + \frac{1}{\varepsilon^2}\|w_m(t)\|^2) + \frac{1}{2\varepsilon^2}\|v_1(t)\|^2_{H^{-1/2}(\Gamma_1)} \tag{3.7}$$

Setting $\varepsilon^2 = c_1/(2\beta\|g\|_\infty + c_{10})$ in (3.7) gives

$$\frac{d}{dt}|z_m(t)|^2 + \frac{c_1}{2}\|z_m(t)\|^2$$

$$\leq \frac{2(\beta\|g\|_\infty + c_{10})}{c_1} \beta\|g\|_\infty \|w_m(t)\|^2 + \frac{2(\beta\|g\|_\infty + c_{10})}{c_1} \|v_1(t)\|^2_{H^{-1/2}(\Gamma_1)} \tag{3.8}$$

By assumption (2.6) again, for any given $\varepsilon > 0$, we can get from (3.5) that

$$\frac{d}{dt}|w_m(t)|^2 + c_1'\|w_m(t)\|^2 \leq \varepsilon^2 \|w_{mn}(t)\|^2_{H^{1/2}(\Gamma_2)} + \frac{1}{\varepsilon^2} \|v_2(t)\|^2_{H^{-1/2}(\Gamma_2)} \tag{3.9}$$

By the trace theorem from $H^1(\Omega)$ to $H^{1/2}(\Gamma)$ again,, there exists a constant $c_{11}$ such that

$$\|w_{mn}(t)\|_{H^{1/2}(\Gamma_2)} \leq c_{11}\|w_m(t)\|$$

Substitute above into (3.9) to yield

$$\frac{d}{dt}|w_m(t)|^2 + (c_1' - c_{11}\varepsilon^2)\|w_m(t)\|^2 \leq \frac{1}{\varepsilon^2}\|v_2(t)\|^2_{H^{-1/2}(\Gamma_2)} \tag{3.10}$$

Setting $\varepsilon^2 = c_1'/(2c_{11})$ in (3.10) gives

$$\frac{d}{dt}|w_m(t)|^2 + \frac{c_1'}{2}\|w_m(t)\|^2 \leq \frac{2c_{11}}{c_1'}\|v_2(t)\|^2_{H^{-1/2}(\Gamma_2)} \tag{3.11}$$

Integrate (3.11) over [0; T] with respect to t to give

$$|w_m(T)|^2 + \frac{c_1'}{2}\int_0^T \|w_m(t)\|^2 dt \leq \frac{2c_{11}}{c_1'}\int_0^T \|v_2(t)\|^2_{H^{-1/2}(\Gamma_2)} dt + |w_m(0)|^2$$

$$\leq \frac{2c_{11}}{c_1'}\int_0^T \|v_2(t)\|^2_{H^{-1/2}(\Gamma_2)} dt + |w(0)|^2 \tag{3.12}$$

Since the right-hand side of (3.12) is bounded, we have

$$\{w_m\} \text{ is a bounded sequence in } L^2(0,T;W). \tag{3.13}$$

Replace T by $t \in [0,T]$ in (3.12) to obtain

$$ess\sup_t |w_m(t)|^2 \leq \frac{2c_{11}}{c_1'}\int_0^T \|v_2(t)\|^2_{H^{-1/2}(\Gamma_2)} dt + |w_m(0)|^2 \tag{3.14}$$

Hence

$$\{w_m\} \text{ is a bounded sequence in } L^\infty(0,T;\widetilde{H}). \tag{3.15}$$

On the other hand, integrate (3.8) over [0; T] with respect to $t$ to give

$$|z_m(T)|^2 + \frac{c_1}{c}\int_0^T \|z_m(t)\|^2 dt \leq \frac{2(\beta\|g\|_\infty + c_{10})}{c_1}\beta\|g\|_\infty \int_0^T \|w_m(t)\|^2 dt +$$

$$\frac{2(\beta\|g\|_\infty + c_{10})}{c_1}\int_0^T \|v_1(t)\|^2_{H^{-1/2}(\Gamma_1)} dt + |z(0)|^2 \tag{3.16}$$

Therefore

$$\{z_m\} \text{ is a bounded sequence in } L^2(0,T;V). \tag{3.17}$$

Replace T by $t \in [0,T]$ in (3.16) to get

$$ess\sup_t |w_m(t)|^2 \leq \frac{2(\beta\|g\|_\infty + c_{10})}{c_1}\beta\|g\|_\infty \int_0^T \|w_m(t)\|^2 dt +$$

$$\frac{2(\beta\|g\|_\infty + c_{10})}{c_1}\int_0^T \|v_1(t)\|^2_{H^{-1/2}(\Gamma_1)} dt + |z(0)|^2 \tag{3.18}$$

Therefore,

$$\{z_m\} \text{ is a bounded sequence in } L^\infty(0,T;H). \tag{3.19}$$

(3.3) and (3.4) then follow from (3.13), (3.17) (3.15) and (3.19).

**Lemma 3.2.** Let $\{z_m, w_m\}$ be the sequence determined by Lemma 3.1. Then there exists a sequence of $\{z_m, w_m\}$, still denoted by itself without confusion, such that

$$z_m \to z \text{ strongly in } L^2(0,T;H), \quad w_m \to w \text{ strongly in } L^2(0,T;\tilde{H}), \tag{3.20}$$

Proof. By virtue of Lemmas 2.1-2.3, we can write (3.2) as follows:

$$\begin{cases} (\dfrac{dz_m(t)}{dt}, u_j) = (H_1 - \beta g w_m(t) - \nu A_1 z_m(t) - B(z_m(t)), u_j), \forall j = 1, 2, \cdots, m, \\ (\dfrac{dw_m(t)}{dt}, \mu_j) = (H_2 - C(z_m,(t) w_m(t)) - k A_2 w_m(t), \mu_j), \forall j = 1, 2, \cdots, m, \end{cases} \tag{3.21}$$

Denote by $\{\tilde{z}_m, \tilde{w}_m\}$ the $\{z_m, w_m\}$ with zero values outside of $[0; T]$ and $\{\hat{z}_m, \hat{w}_m\}$ the Fourier transformations of $\{\tilde{z}_m, \tilde{w}_m\}$. We claim that there exists a $\nu > 0$ such that

$$\int_{-\infty}^{+\infty} |\tau|^{2\nu} \|\hat{z}_m(\tau)\|^2 d\tau < \infty. \tag{3.22}$$

To this end, we write the _rst equation of (3.21) as

$$\frac{d}{dt}(\tilde{z}_m, u_j) = <\tilde{f}_{m1}, u_j> + (z_{0m}, u_j)\delta_0 - (z_m(T), u_j)\delta_T \tag{3.23}$$

where $\delta_0, \delta_T$ are Dirac functions, and

$$\tilde{f}_{m1}(t) = f_{m1}(t) \text{ for } t \in [0,T] \text{ and } \tilde{f}_{m1}(t) = 0 \text{ for } t > T$$
$$f_{1m}(t) = H_1 - \beta g w_m(t) - \nu A_1 z_m(t) - B(z_m(t)),$$

Take Fourier transform for Equation (3.23) to get

$$2\pi i \tau(\tilde{z}_m, u_j) = <\hat{f}_{m1}, u_j> + (z_{0m}, u_j)\delta_0 - (z_m(T), u_j) e^{-2\pi i T\tau} \tag{3.24}$$

where $\hat{f}_{m1}$ is the Fourier transform of $\tilde{f}_{m1}$.

Let $\tilde{q}_{jm}(t)$ be the function of $q_{jm}(t)$ in (3.1), which is zero outside of $[0; T]$ and let $\hat{q}_{jm}(t)$ be its Fourier transform. Multiplier Equation (3.24) by $\hat{q}_{jm}(t)$ and sum for $j$ from 1 to m to obtain

$$2\pi i |\hat{z}_m(\tau)|^2 = <\hat{f}_{m1}, \hat{z}_m(\tau)> + (z_{0m}, \hat{z}_m(\tau))\delta_0 - (z_m(T), \hat{z}_m(\tau)) e^{-2\pi i T\tau} \tag{3.25}$$

We thus conclude that

$$\int_0^T \|f_{1m}(t)\|_{V^*}^2 dt \leq \int_0^T [\|f_1(t)\|_{V^*} + c_2\|w_m(t)\| + c_1\|z_m(t)\| + c_1\|z_m(t)\|^2] dt \tag{3.26}$$

where we used the fact $\|B(z_m(t)\| \leq c_1\|z_m(t)\|^2$. By (3.13), (3.15) and (3.17), it follows from (3.26) that

$$\sup_{\tau \in R} \|\hat{f}_{1m}(t)\|_{V^*} < \infty \tag{3.27}$$

Apply (3.27) and the fact $\sup_{m \in Z^+}[\|z_m(0)\| + \|z_m(T)\|] < \infty$ to (3.25) to yield

$$\tau |\hat{z}_m(\tau)|^2 \leq c_3' \|\hat{z}_m(\tau)\| + c_4' |\hat{z}_m(\tau)| = c_5' \|\hat{z}_m(\tau)\|, c_5' = c_3' + c_4' \tag{3.28}$$

For fixed $0 < \nu < 1/4$, observe that

$$|\tau|^{2\nu} \leq c'(\nu) \frac{1+|\tau|}{1+|\tau|^{1-2\nu}}, \forall \tau \in R$$

for some constant $c'(\nu)$. From this inequality, we obtain

$$\int_{-\infty}^{+\infty} |\tau|^{2\nu} |\hat{z}_m(\tau)|^2 d\tau \leq c'(\nu) \int_{-\infty}^{+\infty} \frac{1+|\tau|}{1+|\tau|^{1-2\nu}} |\hat{z}_m(\tau)|^2 d\tau \tag{3.29}$$

By (3.28), there are constants $c_0' > 0$ and $c_1' > 0$ such that

$$\int_{-\infty}^{+\infty} |\tau|^{2\nu} |\hat{z}_m(\tau)|^2 d\tau \leq c'(\nu) \int_{-\infty}^{+\infty} \frac{1+|\tau|}{1+|\tau|^{1-2\nu}} |\hat{z}_m(\tau)|^2 d\tau$$

$$\leq c_0' \int_{-\infty}^{+\infty} \frac{\|\hat{z}_m(\tau)\|}{1+|\tau|^{1-2\nu}} d\tau + c_1' \int_{-\infty}^{+\infty} \|\hat{z}_m(\tau)\|^2 d\tau$$

By (3.20) and the Parseval identity, the second term on the right-hand side of above inequality is bounded as $m \to \infty$. Therefore, (3.22) is proved if we can show that

$$\int_{-\infty}^{+\infty} \frac{\|\hat{z}_m(\tau)\|}{1+|\tau|^{1-2\nu}} d\tau < \infty$$

However, this is a consequence of the Schwarz inequality and the Parseval identity that

$$\int_{-\infty}^{+\infty} \frac{\|\hat{z}_m(\tau)\|}{1+|\tau|^{1-2\nu}} d\tau < \left(\int_{-\infty}^{+\infty} \frac{d\tau}{1+|\tau|^{1-2\nu}}\right)^{1/2} \left(\int_{-\infty}^{+\infty} \|\hat{z}_m(\tau)\|^2 d\tau\right)^{1/2},$$

where we used the facts $0 < \nu < 1/4$ and the boundedness of $\{z_m\}$ in $L^\infty(0,T;H)$ as $m \to \infty$ claimed by (3.19). So (3.22) is valid.

By (3.22), (3.17) and (3.19), we conclude that

$$\{z_m\} \text{ is bounded in } H^\nu(R;V) \cap H^\nu(R;H) \tag{3.30}$$

By Lemma 2.5, there exists a subsequence of $\{z_m, w_m\}$ that is still denoted by itself without confusion such that

$$z_m \to z \text{ strongly in } L^2(0,T;H), \quad w_m \to w \text{ strongly in } L^2(0,T;\tilde{H}).$$

This is (3.20).

**Theorem 3.1.** There exists a weak solution to (1.1).

Proof. Let $\Psi$ and $\theta$ and be continuous differentiable vector functions defined on $[0; T]$ with $\Psi(T) = \theta(T) = 0$. Multiply the first equation of (3.2) by $\Psi$ and integrate over $[0; T]$ with respect to $t$ to give

$$-\int_0^T (z_m(t), \Psi'(t)u_j) dt + \int_0^T [\nu a_1(z_m(t), \Psi(t)u_j) + b(z_m(t), z_m(t), \Psi(t)u_j) + \beta(w_m(t)g, \Psi(t)u_j)] dt$$

$$= (z_{0m}, u_j)\Psi(0) + \int_0^T <v_1(t), \Psi(t)u_{jn}>_{\Gamma_1} dt \tag{3.31}$$

Multiply the second equation of (3.2) by $\theta$ and integrate over $[0; T]$ with respect to $t$ to give

$$-\int_0^T (w_m(t), \theta'(t)\mu_j) dt + \int_0^T [k a_2(w_m(t), \theta'(t)\mu_j) + c(z_m(t), w_m(t), \theta(t)\mu_j)] dt$$

$$= (w_{0m}, \mu_j)\theta(0) + \int_0^T <v_2(t), \theta(t)\mu_j>_{\Gamma_2} dt \tag{3.32}$$

Passing to the limit as $m \to \infty$ in (3.31) and (3.32) by applying (3.3), (3.20), the properties (iii) in Lemmas 2.2 and 2.3 for $b$ and $c$, we obtain

$$-\int_0^T (z(t), \Psi'(t)u_j) dt + \int_0^T [\nu a_1(z(t), \Psi(t)u_j) + b(z(t), z(t), \Psi(t)u_j) + \beta(w(t)g, \Psi(t)u_j)] dt$$

$$= (z_0, u_j)\Psi(0) + \int_0^T <v_1(t), \Psi(t)u_{jn}>_{\Gamma_1} dt \tag{3.33}$$

$$-\int_0^T (w(t), \theta'(t)\mu_j)dt + \int_0^T [ka_2(w(t), \theta'(t)\mu_j) + c(z(t), w(t), \theta(t)\mu_j)]dt$$

$$= (w_0, \mu_j)\theta(0) + \int_0^T <v_2(t), \theta(t)\mu_j>_{\Gamma_2} dt \qquad (3.34)$$

where in obtaining (3.33) and (3.34), we used the following facts:
- The convergence of the nonlinear terms in $b(\cdot,\cdot,\cdot)$ and $c(\cdot,\cdot,\cdot)$ can be obtained in the same way as that in Chapter 3 of [12].
-
$$\int_0^T a_1(z_m(t), \Psi(t)u_j)dt = \int_0^T (rotz_m(t), \Psi(t)rotu_j)dt =$$

$$\to \int_0^T (rotz(t), \Psi(t)rotu_j)dt = \int_0^T a_1(z(t), \Psi(t)u_j)dt$$

where we used the facts that $rotz_m(t) \to rotz(t)$ weakly in $L^2(0,T;V)$.
- Similarly

$$\int_0^T a_2(w_m(t), \theta(t)\mu_j)dt \to \int_0^T a_2(w(t), \theta(t)\mu_j)dt$$

by the facts again $\nabla w_m(t) \to \nabla w(t)$ weakly in $L^2(0,T;W)$

By the density arguments, we have that (3.33) and (3.34) hold true for any $\psi \in V$ instead of $u_j$ and $\varphi \in W$ instead on $\mu_j$, respectively. That is,

$$-\int_0^T (z(t), \Psi'(t)\psi)dt + \int_0^T [va_1(z(t), \Psi(t)\psi) + b(z(t), z(t), \Psi(t)\psi) + \beta(w(t)g, \Psi(t)\psi)]dt$$

$$= (z_0, \psi)\Psi(0) + \int_0^T <v_1(t), \Psi(t)\psi_n>_{\Gamma_1} dt, \forall \psi \in V \qquad (3.35)$$

$$-\int_0^T (w(t), \theta'(t)\varphi)dt + \int_0^T [ka_2(w(t), \theta'(t)\varphi) + c(z(t), w(t), \theta(t)\varphi)]dt$$

$$= (w_0, \varphi)\theta(0) + \int_0^T <v_2(t), \theta(t)\varphi>_{\Gamma_2} dt \qquad (3.36)$$

Now take $\Psi \in (D(0,T))^N$ in (3.35) and $\theta \in D(0,T)$ in (3.36). Then $\{z, w\}$ satisfies

$$\begin{cases} (z', \psi) + v\, a_1(z, \psi) + b(z, z, \psi) + (\beta g\, w, \psi) = <v_1, \psi_n>_{\Gamma_1}, & \forall \psi \in V \\ (w', \varphi) + ka_2(w, \varphi) + c(z, w, \varphi) = <v_2, \varphi>_{\Gamma_2}, & \forall \varphi \in W \end{cases} \qquad (3.37)$$

This is the equations in (2.1). Finally, we determine the initial value of $\{z, w\}$. Actually, multiply the _rst equation of (3.37) and integrate over [0; T] with respect to $t$ to get

$$-\int_0^T (z(t), \Psi'(t)\psi)dt + \int_0^T [va_1(z(t), \Psi(t)\psi) + b(z(t), z(t), \Psi(t)\psi) + \beta(w(t)g, \Psi(t)\psi)]dt$$

$$= (z(0), \psi)\Psi(0) + \int_0^T <v_1(t), \Psi(t)\psi_n>_{\Gamma_1} dt, \forall \psi \in V$$

Subtract (3.38) from (3.35) to get $(z(0) - z_0, \psi)\Psi(0) = 0$.. Take $\Psi$ so that $\Psi(0) = 1$ to get $(z(0) - z_0, \psi) = 0$ for all $\forall \psi \in V$. So $z(0) = z_0$. The similar arguments lead to $w(0) = w_0$ The proof is complete. $\square$

# 4. Existence of the optimal pair

**[Definition 2]** Suppose that $z_0 \in H$, $w_0 \in \widetilde{H}$, $g \in L^\infty(\Omega)$.
The pair $\{y, v\} = \{(z, w), (v_1, v_2)\}$ is said to be admissible pair; "state-control" of extreme value problem (1.2) if it satisfies

$$\begin{cases} y = \{z, w\} \in Y, z' \in L^1(0, T, V^*), w' \in L^1(0, T : W^*) \\ \dfrac{dz}{dt} + \nu A_1 z + B(z) + \beta g w = H_1 v_1, \\ \dfrac{dw}{dt} + k A_2 w + C(z, w) = H v_2, \\ z(0) = z_0, w(0) = w_0 \\ v = \{v_1, v_2\} \in U_a \equiv U_{1a} \times U_{2a} \end{cases} \quad (4.1)$$

We denote admissible pair sets by M that is;

$$M \equiv \{(y, v) \mid y \equiv \{z, w\} \in Y = Z \times W, v = \{v_1, v_2\} \in U_a, (y, v) \text{ satisfy } (4.1)\}.$$

**[Definition 3]** The admissible pair $\{y_*, v_*\} = \{\{z_*, w_*\}, \{v_1^*, v_2^*\}\} \in M$ is called an optimal pair if $v_*$ is satisfied (1.2) and it is called by optimal control and $y_*$ is called by optimal state.

The optimal control problem (1.2) is an optimal control of the singular distributed system ([12]).
We can easily see that energy inequalities

$$\frac{1}{2}|z(t)|^2 + \nu c_1 \int_0^T \|z(\tau)\|^2 d\tau + \int_0^T (\beta g w(\tau), z(\tau)) d\tau \leq \frac{1}{2}|z_0| + \int_0^T <H_1 v_1, z(\tau)> d\tau \quad (4.2)$$

$$\frac{1}{2}|w(t)|^2 + k c_1' \int_0^T \|w(\tau)\|^2 d\tau \leq \frac{1}{2}|w_0| + \int_0^T <H_2 v_2, w(\tau)> d\tau \quad (4.3)$$

are holed.
From the inequalities (4.2), (4.3), we can obtain as following inequalities;

$$\|z\|^2_{L^\infty(0,T;H)} + \|z\|^2_{L^2(0,T;V)} + \left\|\frac{dz}{dt}\right\|^2_{L^1(0,T;V^*)} \leq \widetilde{c}_1(\|z_0\|^2 + \|v_1\|^2_{L^2(0,T;\Gamma_1)}) \quad (4.4)$$

$$\|w\|^2_{L^\infty(0,T;\widetilde{H})} + \|w\|^2_{L^2(0,T;W)} + \left\|\frac{dw}{dt}\right\|^2_{L^1(0,T;W^*)} \leq \widetilde{c}_2(\|w_0\|^2 + \|v_2\|^2_{L^2(0,T;\Gamma_2)}) \quad (4.5)$$

where and throughout in the sequel of the proof of theorems, $\widetilde{c}_i > 0$ denote positive constants independence $t, z, w$.

Here and in what follows, we denote $\|v_1\|^2_{L^2(0,T;(L^2(\Gamma_1))^N)}$ by $\|v_1\|^2_{L^2(0,T;\Gamma_1)}$ and $\|v_2\|^2_{L^2(0,T;L^2(\Gamma_1))}$ by $\|v_2\|^2_{L^2(0,T;\Gamma_2)}$ simply.

**[Lemma 4.1.]** The admissible pair set $M$ is nonempty and weakly closed subset in the space $[L^2(\Sigma_1)]^N \times L^2(\Sigma_2) \times Y$.

Proof. From theorem 3.1, we obtain that the admissible pair set $M$ is nonempty subset in the space $[L^2(\Sigma_1)]^N \times L^2(\Sigma_2) \times Y$.
Now, we shall prove that M is a weakly closed set.
We assume that $\{\{z_m, w_m\}, \{v_{1m}, v_{2m}\}\} \subset M$ are satisfied as following;

$$v_{1m} \to v_1, \text{ weakly in } L^2(0, T; (L^2(\Gamma_1))^N) \;;\; v_{2m} \to v_2, \text{ weakly in } L^2(0, T; L^2(\Gamma_2))$$
$$(4.5)$$

$$z_m \to z, \text{ weakly in } L^2(0, T; V) \;;\; w_m \to w, \text{ weakly in } L^2(0, T; W) \quad (4.6)$$

From (4.4) we obtain that $\left\|\dfrac{dz_m}{dt}\right\|_{L^1(0,T;V^*)} \leq const$. Now, by using compactness theorem of embedding of in [13, 14], we will show that we can chose the subsequence $\{z_m\}$, which we denote by the same symbols, such that

$$z_m \to z, \text{ strongly in } L^2(0,T:H) \qquad (4.7)$$

Same as above from (4.5) we obtain that $\left\|\dfrac{dw_m}{dt}\right\|_{L^1(0,T;W^*)} \leq const$. Now, by using compactness theorem of embedding of [13, 14], we will show that we can chose the subsequence $\{w_m\}$, which we denote by the same symbols, such that

$$w_m \to w, \text{ strongly in } L^2(0,T:\tilde{H}) \qquad (4.8)$$

From the (4.5)-(4.8) we can take the limit as $m \to \infty$ in the equations;

$$\frac{dz_m}{dt} + \nu A_1 z_m + B(z_m) + \beta g w_m = H_1 v_{1m},$$

$$\frac{dw_m}{dt} + k A_2 w_m + C(z_m, w_m) = H v_{2m}$$

and we obtain that $\{z, w\}$ are satisfied the equations (4.1).

We can rewrite the inequalities (4.2), (4.3) for $\{z_m, w_m\}$, $\{v_{1m}, v_{2m}\}$, and $v_i \in [\alpha_i(x,t), \beta_i(x,t)]$. Therefore, $\{\{z, w\}, \{v_1, v_2\}\}$ is an admissible pair. Thus, proof of lemma is fulfilled. □

**[Theorem 4.1]** There exists optimal pair of the optimization problem (1.2).

Proof. We assume that $\{\{z_m, w_m\}, \{v_{1m}, v_{2m}\}\} \in M$ is the minimization sequences, that is, when $m \to \infty$

$$I[v_m] = I[v_{1m}, v_{2m}] \to \inf I$$

Then, we have obtain such as estimations

$$\|v_{1m}\|_{L^2(\Sigma_1)} \leq const, \forall m \ ; \ \|v_{2m}\|_{L^2(\Sigma_2)} \leq const, \forall m$$

By estimations (4.4), (4.5), we have

$$\|z_m\|_{L^2(0,T;V)} \leq const, \forall m \ ; \ \|w_m\|_{L^2(0,T;W)} \leq const, \forall m$$

We can chose the subsequence $\{w_m\}$, which we denote by the same symbols, by taking the limit as $m \to \infty$, obtain such that

$$v_{1m} \to v_1, \text{ weakly in } L^2(\Sigma_1); \quad v_{2m} \to v_2, \text{ weakly in } L^2(\Sigma_2)$$
$$z_m \to z, \text{ weakly in } L^2(0,T;V); \quad w_m \to w, \text{ weakly in } L^2(0,T;W)$$

Seeing that the admissible pair set M is the nonempty, weakly compactness and cost functional $I$ is continuous, $\{\{z,w\},\{v_1, v_2\}\}$ is an optimal pair by generalized Weierstrass theorem. Thus, proof of theorem is fulfilled. □